\input{style/arxiv-ba.cfg}
\documentclass[ba,linksfromyear,preprint]{imsart}
\usepackage{vtexurl}
\startlocaldefs

\let\bid@start\@empty
\let\bid@end\@empty

\def\MR@url{http://www.ams.org/mathscinet-getitem?mr=}
\def\MR#1{\href{\MR@url#1}{MR#1}}
\def\BDOI#1{%
\edef\doi@base@i{\doi@base}\def\doi@base{}%
doi:~\doiurl{\doi@base@i#1}}

\appto\bid@start{\def\doi@size{\ttfamily}}
\appto\bid@end{\unskip.}

\define@key{bid}{mr}{\def\bid@mr{\MR{#1}}}
\define@key{bid}{doi}{\def\bid@doi{\BDOI{#1}}}
\define@key{bid}{pubmed}{}
\define@key{bid}{pii}{}
\define@key{bid}{issn}{}

\def\bid#1{%
       \bgroup
       \bid@start
       \let\bid@output\@empty
       \setkeys{bid}{#1}\ignorespaces%
       \ifdefvoid\bid@mr{}{\appto\bid@output{\bid@mr}}%
       \ifdefvoid\bid@doi{}{
         \ifdefempty\bid@output{}{\appto\bid@output{. }}% separator if output not empty
         \appto\bid@output{\bid@doi}%
       }%
       \bid@output
       \bid@end
       \egroup
}
\endlocaldefs

\pubyear{2015}
\volume{10}
\issue{1}
\firstpage{243}
\lastpage{246}
\doi{10.1214/15-BA943}% Updated by VTEXPTS2LaTeX.exe, 02.02.2015 16:08

\begin{document}

\begin{frontmatter}
\title{Rejoinder\thanksref{T1}}
\runtitle{Rejoinder}

\relateddois{T1}{Main article DOI: \relateddoi[ms=BA915]{Related item:}{10.1214/14-BA915}.}

\begin{aug}
\author[a]{\fnms{James O.} \snm{Berger}\ead[label=e1]{berger@stat.duke.edu}},
\author[b]{\fnms{Jose M.} \snm{Bernardo}\corref{}\ead[label=e2]{jose.m.bernardo@uv.es}},
\and
\author[c]{\fnms{Dongchu} \snm{Sun}\ead[label=e3]{sund@missouri.edu}}

\runauthor{J. O. Berger, J. M. Bernardo, and D. Sun}

\address[a]{Duke University, USA and King Abdulaziz University, Saudi Arabia, \printead{e1}}
\address[b]{Universitat de Val{\`e}ncia, Spain, \printead{e2}}
\address[c]{University of Missouri-Columbia, USA, \printead{e3}}

\end{aug}

\end{frontmatter}

%% Mainmatter %%

Our thanks to all the discussants for their insightful observations and
comments. We respond to their discussions in turn.

\section{Response to Datta and Liseo}

We agree that Method 3 is preferable to Method 2, in that it is not
dependent on the specification of a collection of quantities of
interest and, hence, need only be determined once (and not separately
for each potential user of a model). It is because a hierarchical
embedding is not always available that we introduce the other methods
as possible solutions.

We found the discussion of the multinomial example interesting, with
numerous additional insights being provided. Likewise the additional
material on the geometric averaging approach was enlightening,
especially the nice lemma showing that, if a collection of priors all
yield proper posteriors, then their geometric average also yields a
proper posterior. This certainly strengthens the argument that
geometric averaging is superior to arithmetic averaging in the search
for an overall prior.

The moral of the amusing anecdote is indeed sound, and can be attempted
to be implemented even when there is no hierarchical embedding
available. For instance, %Berger and Sun (2008)
\cite{BerSun2008} considered 21 different
derived parameters for the five-parameter bivariate normal
distribution, seeking a prior that was good `on average' for the 21
parameters.

\section{Response to Mendoza and Guti{\'e}rrez-Pe{\~n}a}

The discussants highlight the importance of cataloguing those
situations in which there is a common reference prior for all the
parameters of a model and give useful references that could be a
starting point for identifying additional such situations. But they
then, interestingly, question whether this is sufficient, especially
when the number, $m$, of quantities of interest exceeds the number,
$k$, of parameters in the model.

Section 3.1 highlights one such situation: there is a common reference
prior for $\mu$ and $\sigma$ from the normal model but this cannot
necessarily be claimed to be the overall objective prior because the
reference prior for $\mu/\sigma$ is different. This simple example
suggests that one can probably never have an overall objective prior
that is optimal for everything and that just having it be reasonable
for everything (of interest) might be the best we can hope for. It is
interesting, in this regard, that Section 3.2.4 shows that the common
reference prior for $\mu$ and $\sigma$ is also nearly optimal when
$\mu/\sigma$ is also considered, although the discussants are correct
that this result was probably inevitable once we restricted the
candidate priors to be only of the form $\sigma^{-a}$. Utilizing the
alternative and more general class that they suggest might well have
given a different result (but the computation would have been much more
formidable).

It would be nice if one could show, in general, that, if there is a
common reference prior for all of original parameters, then that prior
will be reasonable for other derived parameters or quantities. Our
experience strongly supports this claim, but it is difficult to see how
to formally approach verification of the claim especially, as the
discussants note, because of the disquieting result in Section 2.2.2.

We agree that it would be nice if the family of candidate priors
considered in both the reference distance method and the hierarchical
method could somehow be intrinsically identified from the model itself;
this would make the label `objective' more compelling. We have not
tried to do so ourselves, but the discussants give several potentially
useful starting points for such an endeavor. Computational
considerations are central here so, as the discussants note, we always
chose the candidate class of priors to be a conjugate class (or as
close to conjugate as possible).

Thanks for pointing out the possible relationship of the reference
distance method with the mean field approach to variational inference.
The approximation tools being used in each case are clearly related,
but it is not clear to us that this can be usefully exploited.

We agree with the discussants concerning Section 4.2. It is hard to
know how to deal with the hypergeometric parameters directly, so we
used the common technique of `transferring' them into uncertain
multinomial parameters that we can deal with. But this is, indeed, a
somewhat {\it ad hoc} addition to the proposed methodology. In this
light, the suggested reformulation of the discussants (which ends up in
the same place) will be a more appealing justification to many.

\section{Response to Rousseau}

Rousseau makes the important observation that we are considering the
`simple' parametric case, where there is some hope of having an overall
objective prior that is at least reasonable for likely quantities of
interest. This hope could well be impossible in nonparametric
situations, where it can be a challenge to even find a prior that is
satisfactory for a single given quantity of interest.

Rousseau observes that maybe the search for an overall prior should be
considered together with choice of the model. This is an intriguing
idea, but we have no idea how to approach the issue.

Rousseau observes that, for the reference distance method, the solution
depends on the sample size. It is not appealing, in general, to have
objective priors depend on the sample size, but there are situations
(hierarchical models) where it seems correct and inevitable. Here,
however, the numerical evidence in the examples indicates that there is
only a very slight dependence on the sample size, so Rousseau observes
that one can simply try to implement the approach asymptotically,
avoiding the sample size dependence and -- more importantly -- perhaps
considerably simplifying the derivation of the overall prior. This is
an idea definitely worth pursuing!

In her final comments, Rousseau addresses scenarios considerably more
complex than any we consider, and outlines issues in finding good
(objective) priors for those scenarios. In our own statistical
practice, we encounter these problems all the time. There is little or
no theory to guide us, so it is perhaps most useful to simply say what
we do. A complex model is usually made up of simpler subcomponents, and
we may well know a good overall objective prior for a subcomponent. We
will use it, even though there is no assurance that it is a good
overall objective prior in the context of the full model; the
alternative of using
a prior that we know is suboptimal for the component does not appeal.

This is the but the tip of an iceberg, however, in that many complex
models are hierarchical in nature, and it is well known that standard
objective priors for a model can be terrible if that model appears at a
higher level in a hierarchy. See %Berger, Strawderman, and Tang (2005)
\cite{BerStrTan2005} for discussion of this.

\section{Response to Sivaganesan}
We enjoyed Sivaganesan's comment that ``How \ldots[reference
priors]\ldots\
work seems to be a mystery \ldots,'' because it is also a mystery to us.
But their consistently astonishing properties explains why the
approaches we suggest for developing an overall prior all center around
some application of reference prior theory.

Sivaganesan points out that the choice of the candidate priors will
surely affect the answers, and asks if we have tried alternative
classes of candidate priors. He is certainly right that the class will
likely have some effect, but our experience with Bayesian robustness in
other contexts suggests that the class may not be that important when,
as here, we are optimizing over the class. But this an important topic
for future study.

We appreciated Sivaganesan's comment that ``It is surprising that the
reference prior for $a$ in the hierarchical approach to the multinomial
example turns out to be a proper prior, making up for the behavior of
the marginal [likelihood] being bounded away from~0 at infinity.'' The
history was that we first -- and to our surprise -- discovered that the
reference prior for $a$ was proper; we then went back to look at the
likelihood, and discovered that, indeed, it was not integrable at
infinity, `explaining' why the reference prior decided to be proper.
(This is part of the mystery of reference priors referred to
above.)\looseness=-1

\section{Closing comments}

Mendoza and Guti{\'e}rrez-Pe{\~n}a comment about the paper ``\ldots\ it
offers more of a brainstorming than a systematic treatment and a
general solution to the problem [of obtaining an overall objective
prior].'' We couldn't agree more.\vadjust{\goodbreak} We have been working on this for
more years than we care to reveal and finally admitted to ourselves
that we were not going to find the general solution to the problem. So
the paper is simply a reflection of what we encountered in attempting
to find a general solution.

Sivaganesan asks us to comment on which of the three approaches to an
overall objective prior we would recommend. Details are given in the
final section of the paper, but it is useful to highlight the main
points (with the caveat of the comment above):
\begin{itemize}
\item If all (natural) parameters of the model have the same
reference prior, use it as the overall objective prior.
\item If one can find a natural and computationally feasible
hierarchical structure for the model parameters, use that,
along with finding the reference prior for the parameters in
the hierarchical structure.
\item If the above are not implementable,
\begin{itemize}
\item Try the reference distance approach; the suggestion of
Rousseau to do so asymptotically is perhaps the first thing
to try here.
\item Try the geometric average of parameter reference priors,
supported by the results of Liseo and Datta.
\end{itemize}
\end{itemize}

%% Appendix %%
% \appendix
% \section{}\label{}

%% Supplement Material %%
% \begin{supplement}
% \sname{}\label{}
% \stitle{}
% \slink[url]{}
% \stitlepost{\\} %jei doi reikia nukelti i kita eilute
% \sdescription{}
% \end{supplement}

%% References %%
%

% Acknowledgements
% \begin{acknowledgement}
% \end{acknowledgement}

\end{document}